\documentclass[a4paper,11pt,reqno]{amsart}

\usepackage{amssymb,amsmath,amsthm,color}

\usepackage[bookmarksopen=true,bookmarksnumbered=true, bookmarksopenlevel=2, colorlinks=true,linkcolor=mblue,
citecolor=mblue,urlcolor=mblue, pdfstartview=FitH]{hyperref}

\definecolor{mblue}{rgb}{0,0,.8}

\newcommand{\Z}{\mathbb Z}
\newcommand{\Q}{\mathbb Q}
\newcommand{\F}{\mathbb F}

\renewcommand{\P}{\mathbb P}
\newcommand{\OO}{\mathcal O}

\newtheorem{thm}{Theorem}

\begin{document}

\title[Explicit computation of cusp forms and its complexity]{Explicit computation of cusp forms via Hecke action on cohomology and its complexity}

\author{Jonas B. Rasmussen}


\begin{abstract}
In the literature, the standard approach to finding bases of spaces of modular forms is via modular symbols and the homology of modular curves. By using the Eichler-Shimura isomorphism, a work by Wang shows how one can use a cohomological viewpoint to determine bases of spaces of cusp forms on $\Gamma_0(N)$ of weight $k \geq 2$ and character $\chi$. It is interesting to look at the complexity of this alternative approach, and we do this for an explicit implementation of the algorithm suggested by Wang.
\end{abstract}

\maketitle

\section{Introduction}
When working in the field of modular forms it is often extremely useful to be able to work explicitly with spaces of modular forms, for instance in working with elliptic curves, testing conjectures, etc.

To the author's knowledge, every available software package uses modular symbols to compute bases of spaces of modular forms, and in the literature this is the standard approach as well. A good reference for the computational aspects of this is Stein \cite{Stein}.

In this paper we look at another way of determining bases of spaces of cusp forms, by using a cohomological approach, based on the Eichler-Shimura isomorphism, as done in Wang \cite{Wang}. The theory behind this algorithm is worked through in \cite{Wang}, and what we do here is take an explicit implementation of this algorithm and analyze its complexity.

Our main result is the following:
\begin{thm} \label{thm}
An upper bound on the theoretical complexity of determining a basis for $S_k(\Gamma_0(N),\chi)$ via the cohomological approach described below is
\[
	\OO\big(N^{3+\epsilon}k^{2+\epsilon}(N+k^4)\big),
\]
for $\epsilon > 0$.
\end{thm}

Finally we give two examples, where we work through the main steps of the algorithm. The first of these examples is the easier case of trivial character, while the second example with non-trivial character showcases some other aspects of the algorithm.

\subsection{Notation} \label{notation}
Let $\Delta$ denote all $2 \times 2$ matrices with coefficients in $\Z$, and let $\Gamma = \operatorname{SL}_2(\Z)$ be the matrices herein with determinant $1$. We also define
\[
\begin{split}
	\Delta_0(N) & = \left\{\begin{pmatrix} a & b \\ c & d \end{pmatrix} \in \Delta \ \bigg| \ c \equiv 0 \bmod N\right\}, \\
	\Gamma_0(N) & = \left\{\begin{pmatrix} a & b \\ c & d \end{pmatrix} \in \Gamma \ \bigg| \ c \equiv 0 \bmod N\right\}. \\
\end{split}
\]

We wish to determine a basis for the space of cusp forms of weight $k \geq 2$ on $\Gamma_0(N)$ with character $\chi$ (a Dirichlet character on $(\Z/N\Z)^*$). We denote by $\mu$ the index of $\Gamma_0(N)$ in $\Gamma$.

For $\delta_0 = \big(\begin{smallmatrix} a & b \\ c & d \end{smallmatrix}\big) \in \Delta_0$, we define $\chi(\delta_0) = \chi(d)$.

As in \cite{Wang}, we let $\P_{\chi}^1(\Z/N\Z)$ denote
\[
	\{(c,d) \in \Z/N\Z \times \Z/N\Z \mid \gcd(c,d,N) = 1\}
\]
modulo the relation
\[
	(\lambda c,\lambda d) \sim \chi(\lambda)(c,d), \ \lambda \in (\Z/N\Z)^*,
\]
and we see that $\P_{\chi}^1(\Z/N\Z)$ is just $\P^1(\Z/N\Z)$ if $\chi$ is the trivial character.

For $\gamma = \big(\begin{smallmatrix} a & b \\ c & d \end{smallmatrix}\big) \in \Gamma$ we put
\[
	\overline{\gamma} = (c \bmod N,d \bmod N) \in \P_{\chi}^1(\Z/N\Z),
\]
and it is easily checked that $\overline{\gamma_0\gamma} = \chi(\gamma_0)\overline{\gamma}$ for $\gamma_0 \in \Gamma_0(N)$ and $\gamma \in \Gamma$.

We have an operation of $\Delta$ on
\[
	M = \left\{\sum_{j=0}^{k-2}a_jx^jy^{k-2-j} \ \bigg| \ a_0,\ldots,a_{k-2} \in \Z\right\}
\]
given by
\[
	\begin{pmatrix} a & b \\ c & d \end{pmatrix}x^jy^{k-2-j} = (ax+by)^j(cx+dy)^{k-2-j}, \ \begin{pmatrix} a & b \\ c & d \end{pmatrix} \in \Delta,
\]
and we get an action of $\Delta_0(N)$ on $M_{\chi}$ (where the coefficient ring of $M$ is extended with the values of $\chi$) by setting $\delta_0.m = \chi(\delta_0)(\delta_0m)$.

We also look at the coinduced module $W_{\chi}$ of $M_{\chi}$ on $\Gamma$:
\[
	W_{\chi} = \big\{w : \Gamma \to M_{\chi} \mid w(\gamma_0\gamma) = \gamma_0.w(\gamma) \text{ for } \gamma_0 \in \Gamma_0(N)\big\}.
\]
We get an action of $\delta \in \Delta$ on $w \in W_{\chi}$ by setting
\[
	(\delta.w)(\gamma) = \begin{cases} 0, & \gamma\delta \notin \Delta_0(N)\Gamma, \\ \delta_0.w(\gamma'), & \gamma\delta = \delta_0\gamma', \delta_0 \in \Delta_0(N),\gamma' \in \Gamma, \end{cases}
\]
for $\gamma \in \Gamma$.

For a matrix $\delta \in \Delta$, we denote by $W_{\chi}^\delta$ the submodule of $W_{\chi}$ invariant under the action of $\delta$.

We reserve $I$, $S$, $Q$ and $\varepsilon$ for the following matrices:
\[
	I = \begin{pmatrix} 1 & 0 \\ 0 & 1 \end{pmatrix}, \ S = \begin{pmatrix} 0 & -1 \\ 1 & 0 \end{pmatrix}, \ Q = \begin{pmatrix} 0 & -1 \\ 1 & 1 \end{pmatrix}, \  \varepsilon = \begin{pmatrix} -1 & 0 \\ 0 & 1 \end{pmatrix}.
\]

\section{Algorithm and implementation}

\subsection{Coset representatives}
We start out by getting coset representatives for $\Gamma_0(N)$ in $\Gamma$. Besides the $N$ representatives $\big(\begin{smallmatrix} 0 & -1 \\ 1 & 0 \end{smallmatrix}\big),\ldots,\big(\begin{smallmatrix} 0 & -1 \\ 1 & N-1 \end{smallmatrix}\big)$, we also have $\big(\begin{smallmatrix} 1 & 0 \\ 0 & 1 \end{smallmatrix}\big)$ (if $N$ is prime these are all of them).

As shown in Cremona \cite{Cremona}, there is a bijection between coset representatives $\big(\begin{smallmatrix} a & b \\ c & d \end{smallmatrix}\big)$ of $\Gamma_0(N)$ in $\Gamma$ and elements $(c,d) \in \P^1(\Z/N\Z)$. To get the remaining representatives we simply take the remaining elements $(c,d)$ of $\P^1(\Z/N\Z)$ and lift these to a matrix $\big(\begin{smallmatrix} a & b \\ c & d \end{smallmatrix}\big)$ in $\Gamma$ via the Euclidean algorithm.

It is described in Section 2.2 of \cite{Cremona} how to efficiently determine $\P^1(\Z/N\Z)$. One starts out with the obvious elements $(1,0)$, $(1,0),\ldots(1,N-1)$, and then look at elements $(c,d)$, where $c \mid N$ and $d=1,\ldots,N-1$, and add it to the list if it is not equivalent to an element already on the list (one uses that two elements $(c_1,d_1)$ and $(c_2,d_2)$ are equivalent if and only if $c_1d_2 \equiv c_2d_1 \pmod N$). See also Section 8.7 of \cite{Stein}.

We will denote the coset representatives by $\gamma_1,\ldots,\gamma_{\mu}$.

\subsection{\texorpdfstring{Action of $\Delta$ on $W_{\chi}$ and relations matrix}{Action of Delta on W and relations matrix}}
From p.\ 103 of \cite{Wang} we have the exact sequence (of complex vector spaces)
\[
	0 \longrightarrow S_k(\Gamma_0(N),\chi) \longrightarrow H^1(\Gamma,W_{\chi})_+ \longrightarrow H^1(\langle T \rangle,W_{\chi})_+ \longrightarrow 0,
\]
where $T = \big(\begin{smallmatrix} 1 & 1 \\ 0 & 1 \end{smallmatrix}\big)$.

By p.\ 105 of \cite{Wang} we have that $H^1(\Gamma,W_{\chi})_+ \cong W_{\chi}/(W_{\chi}^S+W_{\chi}^Q+W_{\chi}^\varepsilon)$, which is $W_{\chi}$ modulo the relations
\[
	w+S.w = w+Q.w+Q^2.w = w+\varepsilon.w = 0.
\]

Therefore we need to be able determine a matrix representation of the action of a matrix in $\Delta$ on $W_{\chi}$.

An element of $W_{\chi}$ is determined by its values on the $\mu$ coset representatives, and since $M_{\chi}$ is generated by the $k-1$ homogenous monomials of degree $k-2$, the space $W_{\chi}$ has $\mu(k-1)$ generators.

As a basis for $W_{\chi}$ we thus have the elements
\[
	w_{ij} : \begin{cases} \gamma_r \mapsto x^jy^{k-2-j}, & r = i, \\ 0, & r \neq i, \end{cases}
\]
with $i=1,\ldots,\mu$ and $j=0,\ldots,k-2$.

To determine the action of $\delta \in \Delta$ on $W_{\chi}$, we only need the action of $\delta.w_{ij}$ on the coset representatives. If $\gamma_r\delta \notin \Delta_0(N)\Gamma$ the action is $0$, so we now assume that we can write $\gamma_r\delta = \delta_0\gamma$ for some $\delta_0 \in \Delta_0(N)$ and $\gamma \in \Gamma$. Since we have $\gamma = \gamma_0\gamma_s$ for some $\gamma_0 \in \Gamma_0(N)$ and a coset representative $\gamma_s$, we replace $\delta_0$ with $\delta_0\gamma_0 \in \Delta_0(N)$, so that the action is given by
\[
	\delta_0.w_{ij}(\gamma_s) = \begin{cases} \delta_0.x^jy^{k-2-j}, & s = i, \\ 0, & s \neq i. \end{cases}
\]

To get the action of $\delta$ on $w_{ij}$ we need to run through the coset representatives $\gamma_r$, get the corresponding $\delta_0 \in \Delta_0(N)$ such that $\gamma_r\delta = \delta_0\gamma_s$, and then compute the coefficients $a_0,\ldots,a_{k-2}$ of the polynomial
\[
	\delta_0.x^jy^{k-2-j} = \sum_{t=0}^{k-2}a_tx^ty^{k-2-t},
\]
for $j=0,\ldots,k-2$. These coefficients are then placed in the $(s+j\mu)$'th row and the $(i+t\mu)$'th columns ($t=0,\ldots,k-2$) of a $\mu(k-1) \times \mu(k-1)$-matrix.

This way we get matrix representations of the actions of $I+S$, $I+Q+Q^2$ and $I+\varepsilon$ on $W_{\chi}$, and $W_{\chi}/(W_{\chi}^S+W_{\chi}^Q+W_{\chi}^\varepsilon)$ is then the nullspace of the resulting relations matrix (the relations matrix is the above three matrix representations stacked on top one another).

\subsection{Dimension}
To determine the dimension of the space of cusp forms, we need the dimension of $H^1(\langle T \rangle,W_{\chi})_+$, that is the dimension of the subspace of $H^1(\langle T \rangle,W_{\chi})$ invariant under the action of $\varepsilon$.

By Lemma 6 of \cite{Wang} we have an isomorphism (of $\Q$-vector spaces)
\[
	H^1(\langle T \rangle,W_{\chi}) \cong \bigoplus_{\text{$s$ cusps}}\Q,
\]
and the action of $\varepsilon$ on $\bigoplus_{\text{$s$ cusps}}\Q$ is shown to be given by
\[
	\varepsilon.\{s\} = -\chi(\delta_0)\{s'\}, \ \gamma_j\varepsilon = \delta_0\gamma_i, \ s = \gamma_i^{-1}\infty, s' = \gamma_j^{-1}\infty,
\]
with $i,j \in \{1,\ldots,\mu\}$.

Let us write $\varepsilon.\{s\} = c\{s'\}$. Since $\varepsilon$ is an involution, we have $\varepsilon.\{s'\} = c^{-1}\{s\}$. If $\{s\} = \{s'\}$, we have $\varepsilon.\{s\} = \pm\{s\}$, showing that $\{s\}$ is in the corresponding $\pm$-space. If $\{s\} \neq \{s'\}$, we have $\{s\}\pm c\{s'\}$ in the corresponding $\pm$-space.

Thus, if $\{s\} \neq \{s'\}$ we get elements of both $\pm$-spaces, but when $\{s\} = \{s'\}$ we get an element of only one of these spaces -- the space corresponding to the sign of $-\chi(\delta_0)$, i.e.\ we get at element of the $+$-space if and only if $\chi(\delta_0) = -1$ when $\{s\} = \{s'\}$.

We therefore need to determine when two cusps $\{\gamma_m^{-1}\infty\}$ and $\{\gamma_n^{-1}\infty\}$ are equivalent, and this happens exactly when $\gamma_mT^u = \gamma_0\gamma_n$ for some $\gamma_0 \in \Gamma_0(N)$ and $u \in \Z$. This groups the coset representatives into $\nu_{\infty}$ classes (one for each cusp equivalence class).

Thus $\varepsilon$ maps a cusp $\{s\}$ to $\pm\{s\}$ (with $s = \gamma_i^{-1}\infty$) if and only if the coset representatives $\gamma_i$ and $\gamma_j$, satisfying $\gamma_j\varepsilon = \delta_0\gamma_i$, are in the same equivalence class, which is exactly when $\gamma_iT^u = \gamma_0\gamma_j$ for some $\gamma_0 \in \Gamma_0(N)$ and $u \in \Z$.

Since we only need one representative for each equivalence class, we choose for a given $\gamma_i$, the unique representative $\gamma_j$ satisfying $\gamma_iT = \gamma_0\gamma_j$ for some $\gamma_0 \in \Gamma_0(N)$, and write $\gamma_j\varepsilon = \delta_0\gamma_i$, checking the sign of $\chi(\delta_0)$.

This way we find the part of the $\pm$-spaces coming from the case where $\varepsilon$ maps $\{s\}$ to $\pm\{s\}$. When this does not happen, we get elements of both $\pm$-spaces, and since the sum of the dimensions is $\nu_{\infty}$ we get the dimension of the $+$-space.

The dimension of $S_k(\Gamma_0(N),\chi)$ is the difference between the dimension of the nullspace of the relations matrix and the dimension just found.

\subsection{Hecke action and basis}
Since $S_k(\Gamma_0(N),\chi)$ is the kernel of the homomorphism
\[
	H^1(\Gamma,W_{\chi})_+ \to H^1(\langle T \rangle,W_{\chi})_+,
\]
we take elements in the kernel of this map and compute the corresponding $q$-expansions until we have enough forms to generate $S_k(\Gamma_0(N),\chi)$.

It is described in \cite{Wang} how to choose elements in the kernel, and we will briefly recount this here.

By Lemmas 3 and 4 in \cite{Wang} we get that $W_{\chi} \cong M \otimes \P_{\chi}^1(\Z/N\Z)$ via the map
\begin{equation} \label{iso}
	w \mapsto \sum_{i=1}^{\mu}\gamma_i^{-1}.w(\gamma_i) \otimes \overline{\gamma}_i. \tag{$\ast$}
\end{equation}

Since $H^1(\langle T \rangle,W_{\chi}) \cong W_{\chi}/(1-T)W_{\chi}$, Lemma 6 of \cite{Wang} gives an isomorphism (of $\Q$-vector spaces) $W_{\chi}/(1-T)W_{\chi} \cong \bigoplus_{\text{$s$ cusps}}\Q$ induced by the map $M \otimes \P_{\chi}^1(\Z/N\Z) \to \bigoplus_{\text{$s$ cusps}}\Q$ given by
\[
	m \otimes (c,d) \mapsto m(0,1)\{\gamma^{-1}.\infty\} = m(0,1)\big\{-\tfrac{d}{c}\big\}, \ \overline{\gamma} = (c,d).
\]

Thus, the homomorphism $H^1(\Gamma,W_{\chi}) \to H^1(\langle T \rangle,W_{\chi})$ becomes a homomorphism $W_{\chi}/(W_{\chi}^S+W_{\chi}^Q) \to \bigoplus_{\text{$s$ cusps}}\Q$, and with the above isomorphisms this map is on $W_{\chi} \cong M \otimes \P_{\chi}^1(\Z/N\Z)$ given by
\begin{equation} \label{map}
	m \otimes (c,d) \mapsto m(0,1)\big\{-\tfrac{d}{c}\big\}-m(1,0)\big\{\tfrac{c}{d}\big\}. \tag{$\dagger$}
\end{equation}

For $m \otimes (c,d)$ to be in the kernel, we can use any $(c,d) \in \P_{\chi}^1(\Z/N\Z)$ with $m$ any (non-empty) linear combination of the monomials $xy^{k-3},\ldots,yx^{k-3}$ if $k>2$, and $m=1$ if $k=2$ (in the weight-$2$ case we have to have $\chi(c) = \chi(d)$ as well).

Since elements $(c,d) \in \P^1(\Z/N\Z)$ correspond bijectively to coset representatives $\big(\begin{smallmatrix} a & b \\ c & d \end{smallmatrix}\big)$ of $\Gamma_0(N)$ in $\Gamma$, we can therefore represent elements of the kernel as $m \otimes \gamma_r$ with $r=1,\ldots,\mu$ and $m$ as above.

We use the Heilbronn-Merel matrices (see Proposition 20 of Merel \cite{Merel})
\[
	H_n = \left\{\begin{pmatrix} a & b \\ c & d \end{pmatrix} \in \Delta \ \bigg| \ ad-bc = n, a > b \geq 0, d > c \geq 0\right\}
\]
to determine the action of the Hecke operator $T_n$ on the elements of the kernel found above.

The action of the Hecke operator $T_n$ on $W_{\chi}/(W_{\chi}^S+W_{\chi}^Q)$ is given by
\[
	T_nx = \sum_{A \in H_n}A.x.
\]

Translating this through the isomorphisms above, the action on a kernel element $m \otimes \gamma_r$ is
\[
	T_n(m \otimes \gamma_r) = \sum_{A \in H_n}\chi(\delta_{0,A})(A.m) \otimes \gamma_{r_A},
\]
where we for each $A \in H_n$ write $\gamma_rA = \delta_{0,A}\gamma_{r_A}$ with $\delta_{0,A} \in \Delta_0(N)$ (terms where it is not possible to write $\gamma_rA$ in this way are ignored).

The coefficients of the polynomial $\chi(\delta_{0,A})(A.m)$ are then saved in a vector $t_n$, where the coefficient of $x^jy^{k-2-j}$ is added to the $r_A(k-1)-(k-2-j)$'th entry, as $A$ runs through $H_n$.

If we wish to compute the basis up to exponent $q^M$, we compute $t_n$ for $n=1,\ldots,M$, where we note that $M$ has to be at least $\lfloor\frac{\mu k}{12}\rfloor$, and we let $t$ be the matrix whose $n$'th column is $t_n$.

We then multiply the nullspace matrix of the relations matrix found earlier with the matrix whose $n$'th column is $t_n$, and denote the resulting matrix by $B$. If $B$ has rank equal to the dimension of $S_k(\Gamma_0(N),\chi)$, we have found a basis (the leading rows of $B$).

If $B$ has rank less than the dimension, we choose another element $m \otimes \gamma_r$ in the kernel, compute the Hecke action $t_1,\ldots,t_M$ on this, multiply the nullspace matrix of the relations matrix with the resulting matrix, and get a matrix whose rows are concatenated to $B$, and we again compute the rank of $B$. This procedure continues until we get as many linearly independent rows in $B$ as the dimension of $S_k(\Gamma_0(N),\chi)$.

Experimentation indicates that if we choose $m = xy^{k-3}+\cdots+x^{k-3}y$ (for $k > 2$) this procedure is likely to give a basis using just the first few $\gamma_r$'s.

\subsubsection{Determining the kernel}
In the implementation described, we choose certain kernel elements and generate Fourier coefficients from these until we have enough forms to generate the space of cusp forms.

However, there is no certainty that this will work, that is there is no guarantee that this approach will give enough linearly independent forms (even though the author has yet to see an example of this).

Another approach (which is certain to work every time) is the following. From the nullspace of the relations matrix we have a basis for $W_{\chi}/(W_{\chi}^S+W_{\chi}^Q+W_{\chi}^\varepsilon)$. By using the isomorphism $W_{\chi} \cong M \otimes \P_{\chi}^1(\Z/N\Z)$ given by (\ref{iso}), we can get the image of this basis on a quotient of $M \otimes \P_{\chi}^1(\Z/N\Z)$, expressed in terms of the standard basis of $M \otimes \P_{\chi}^1(\Z/N\Z)$.

We use this basis to write up a matrix representation of the map (\ref{map}) on this quotient, and determine its nullspace (which has dimension equal to the dimension of $S_k(\Gamma_0(N),\chi)$). From the nullspace matrix we read off a basis, and we then write the kernel elements as linear combinations of the $m \otimes \gamma_r$. We now compute the Hecke action on these elements, which we know will generate enough forms to give a basis for $S_k(\Gamma_0(N),\chi)$.

\section{Complexity of implementation}
We always assume that the level is given via its prime factorization, i.e.\ no work is needed to find the divisors of $N$.

We also use that a Dirichlet character on $(\Z/N\Z)^*$ is defined via a lookup table, which takes $\OO(N)$ to create, but we do not then need to worry about the cost of evaluating the character.

\subsection{Coset representatives}
The number of coset representatives is
\[
	\mu = N\prod_{p \mid N}\big(1+\tfrac{1}{p}\big).
\]

Let $n$ be the number of prime divisors of $N$, and let $p_1,\ldots,p_n$ be the first $n$ primes. By using Landau \cite{Landau}, p.\ 139, we find that
\[
\begin{split}
	\prod_{p \mid N}\big(1+\tfrac{1}{p}\big) & \leq \prod_{p \mid N}\big(1+\tfrac{1}{p}+\tfrac{1}{p^2}+\cdots\big) = \prod_{p \mid N}\big(1-\tfrac{1}{p}\big)^{-1} \leq \prod_{i=1}^n\big(1-\tfrac{1}{p_i}\big)^{-1} \\
	& = \OO(\log p_n) = \OO(\log n) = \OO(\log\log N),
\end{split}
\]
since $N \geq 2^n$. We therefore have $\mu = \OO(N\log\log N)$.

We determine the coset representatives by looking at elements $(c,d)$, where $c \mid N$ and $d=1,\ldots,N-1$. By Theorem 315 of Hardy-Wright \cite{Hardy-Wright}, the number of divisors of $N$ is $\OO(N^\delta)$ for $\delta > 0$, and we therefore look at $\OO(N^{1+\delta})$ elements $(c,d)$.

Everytime we look at an element we check if it is equivalent to something already found, and this takes $\OO(\log^2N)$ each time. Whenever we find a new element we lift it to a matrix in $\Gamma$ via the Euclidean algorithm, which also takes $\OO(\log^2N)$. This gives a total complexity of $\OO(N^{1+\delta}\log^2N)$.

\subsection{\texorpdfstring{Action of $\Delta$ on $W_{\chi}$ and relations matrix}{Action of Delta on W and relations matrix}}

All it takes to get $\gamma_r\delta$ on the form $\delta_0\gamma_s$ is to compute some greatest common divisors and run through the coset representatives to see which one works. All in all the complexity for this is $\OO(\mu\log^2N)$.

The hardest part of computing the action on a matrix on a polynomial is the binomial coefficients that shows up when computing polynomial coefficients. The cost of computing $\binom{n}{m}$ is $\OO(m^2\log^2n)$, and so a rough estimate for computing the action of $\delta_0$ on the $k-1$ monomials $x^jy^{k-2-j}$ is $\OO(k^4\log^2k)$. This needs to be done for every $\gamma_r$, i.e.\ $\mu$ times, giving a total complexity of $\OO(\mu^2\log^2N+\mu k^4\log^2k)$ for determining the matrix representation of $\delta$ on $W_{\chi}$.

In can be noted that if we work over the finite field $\F_p$ instead of the integers, one can use a congruence first proved by Lucas \cite{Lucas} to obtain the complexity $\OO(p^2\log^2p\log k)$ for the binomial coefficient computations. A more modern reference for the Lucas congruence is Stanley \cite{Stanley}, p.\ 44.

The relations matrix consists of the matrix representations of the actions of $I+S$, $I+Q+Q^2$ and $I+\varepsilon$, and so is a matrix of size $3\mu(k-1) \times \mu(k-1)$. To compute the nullspace therefore takes $\OO(\mu^3k^3)$, and this is really the time-consuming function of this part.

\subsection{Dimension}
We first build an array whose $i$'th entry is the index $j$ of the coset representative satisfying $\gamma_iT = \gamma_0\gamma_j$, as well as a similar array giving the index $j$ of the coset representative satisfying $\gamma_i\varepsilon = \delta_0\gamma_j$. Doing this takes $\OO(\mu^2\log^2N)$.

Next we determine to which cusp equivalence class each coset representative belong, and computing $\nu_{\infty}$ along the way. The work needed is already done in the first array we created.

We then choose a representative $\gamma_i$ of a cusp equivalence class and use the second array to find an equivalent representative $\gamma_j$ satisfying $\gamma_j\varepsilon = \delta_0\gamma_i$ for some $\delta_0 \in \Delta_0(N)$. We then compute $\chi(\delta_0)$ and add $1$ to the count of the corresponding $\pm$-variable $d_{\pm}$.

The dimension of the $+$-space is then $d_++(\nu_{\infty}-d_+-d_-)/2$, and the dimension of $S_k(\Gamma_0(N),\chi)$ is the difference between the dimension of the nullspace of the relations matrix and the dimension of the $+$-space.

The work in creating the arrays is by far the most work in this, so the complexity of this algorithm is $\OO(\mu^2\log^2N)$.

\subsection{Hecke action and basis}
Even though we use the Heilbronn-Merel matrices in the implementation, we turn to \cite{Merel} for the complexity analysis, since this paper gives another class of matrices which can be used instead of the Heilbronn-Merel matrices, and we have estimates on the size of these classes.

In Section 3 of \cite{Merel}, Merel defines a set $\mathcal{S}_n$, where a matrix $\big(\begin{smallmatrix} a & b \\ c & d \end{smallmatrix}\big) \in \Delta$ is in $\mathcal{S}_n$ if it has determinant $n$ and at least one of the following conditions are satisfied:
\begin{itemize}
	\item $a > |b|$, $d > |c|$, $bc > 0$,
	\item $b=0$, $|c| < d/2$,
	\item $c=0$, $|b| < a/2$.
\end{itemize}

He also defines a set $\mathcal{S}_n'$, where $\big(\begin{smallmatrix} a & b \\ c & d \end{smallmatrix}\big) \in \Delta$ is in $\mathcal{S}_n'$ if it has determinant $n$ and one of the two following conditions are satisfied:
\begin{itemize}
	\item $b=0$, $|c| = d/2$,
	\item $c=0$, $|b| = a/2$.
\end{itemize}

It is easily seen that an upper bound for $|\mathcal{S}_n'|$ is $2\sigma(n)$, where $\sigma(n)$ is the sum of the positive divisors of $n$, and from p.\ 85 of \cite{Merel} we have, as $n \to \infty$,
\[
	|\mathcal{S}_n| \sim \frac{12\log 2}{\pi^2}\sigma(n)\log n.
\]

By Theorem 322 of \cite{Hardy-Wright} we have $\sigma(n) = \OO(n^{1+\delta})$ for $\delta > 0$, so that we have
\[
	|\mathcal{S}_n \cup \mathcal{S}_n'| = \OO(n^{1+\delta}\log n) = \OO(n^{1+\epsilon}),
\]
for $\epsilon > 0$, since $\log n = \OO(n^{\delta'})$ for $\delta' > 0$.

We need to compute $\mathcal{S}_n$ and $\mathcal{S}_n'$ for $n=1,\ldots,\lfloor\frac{\mu k}{12}\rfloor$ (or more $n$ if we want higher precision), so $\OO(\mu^{2+\epsilon}k^{2+\epsilon})$ matrices are needed to compute the Hecke action.

We want to compute the action of $T_n$ on $m \otimes \gamma_r$, for $m$ a linear combination of the monomials $xy^{k-3},\ldots,x^{k-3}y$ if $k>2$ and $m=1$ if $k=2$. To do this we write, for each $A \in \mathcal{S}_n \cup \mathcal{S}_n'$, $\gamma_rA = \delta_{0,A}\gamma_{r_A}$, with $\delta_{0,A} \in \Delta_0(N)$, and this can be done in $\OO(\mu\log^2N)$ for each $A$. We then have
\[
	T_n(m \otimes \gamma_r) = \sum_{A \in \mathcal{S}_n}\chi(\delta_{0,A})(A.m) \otimes \gamma_{r_A} + \frac{1}{2}\sum_{A \in \mathcal{S}_n'}\chi(\delta_{0,A})(A.m) \otimes \gamma_{r_A},
\]
and as mentioned earlier, the complexity of determining the action of $A$ on a linear combination of all possible monomials is $\OO(k^4\log^2k)$.

The coefficients of $\chi(\delta_{0,A})(A.m)$ are added for each $A$ and saved in a vector $t_n$, with indices depending on the monomial $x^jy^{k-2-j}$ and index $r_A$. Determining all necessary $t_n$ are done in $\OO(\mu^{2+\epsilon}k^{2+\epsilon}(\mu\log^2N+k^4\log^2k))$.

Multiplying the nullspace matrix of the relations matrix with the matrix whose $n$'th column is $t_n$, takes $\OO(\mu^3k^3)$, which is less than the complexity of determining the Hecke action.

The resulting matrix is the basis matrix if it has rank equal to the dimension of $S_k(\Gamma_0(N),\chi)$. We therefore do Gaussian elimination and compute the rank to see if we are done. If not, we choose another coset representative $\gamma_r$, get the resulting $t_n$'s of the Hecke action on $m \otimes \gamma_r$, multiply the nullspace matrix of the relations matrix with the resulting matrix, and get a matrix whose rows are concatenated to $B$. We again do Gaussian elimination and compute the rank, and this is repeated until we get the right rank. If we run through all the coset representatives without getting the right rank, we can try with another $m$.

Gaussian elimination is done in $\OO(\mu^3k^3)$, and is therefore insignificant compared to the computation of the Hecke action.

For a given $m$, this procedure is repeated at most $\mu$ times, but experimentation indicates that it is likely to finish much sooner if $m$ is chosen to be the sum of all possible monomials (in the case of $k > 2$). We see that the computation of the Hecke action is by far the hardest part of this basis determination, and we therefore get a total theoretical complexity of $\OO(\mu^{3+\epsilon}k^{2+\epsilon}(\mu\log^2N+k^4\log^2k))$.

Since $\mu = \OO(N\log\log N)$, and hence $\mu^{\epsilon} = \OO(N^{\epsilon})$, we find that
\[
\begin{split}
	\OO\big(\mu^{3+\epsilon}k^{2+\epsilon}(\mu\log^2N+k^4\log^2k)\big)
	& = \OO\big(\mu^{4+\epsilon}k^{2+\epsilon}\log^2N+\mu^{3+\epsilon}k^{6+\epsilon}\log^2k\big) \\
	& = \OO\big(\mu^{4+\epsilon+\epsilon'}k^{2+\epsilon}+\mu^{3+\epsilon}k^{6+\epsilon+\epsilon'}\big) \\
	& = \OO\big(N^{3+\epsilon+\epsilon'}k^{2+\epsilon+\epsilon'}(N+k^4)\big),
\end{split}
\]
using that $\log^2\alpha = \OO(\alpha^{\epsilon'})$ for $\epsilon' > 0$ (and $\alpha$ either $N$ or $k$), and this is Theorem~\ref{thm} after replacing $\epsilon+\epsilon'$ with $\epsilon$.

\section{Examples}
We use this implementation to determine bases for two spaces of cusp forms, one with trivial character and one with non-trivial character.

The first example gives more detail, while the second highlights an aspect which only happens in the case of non-trivial character.

\subsection{\texorpdfstring{$S_4(\Gamma_0(25))$}{S4(25)}}
We start out by getting the coset representatives, and we do this by determining $\P^1(\Z/25\Z)$ in the way we described earlier. Thus we get the $30$ elements
\[
\begin{split}
	& (1,0),\ldots,(1,24),(0,1),(5,1),(5,2),(5,3),(5,4),
\end{split}
\]
and these are lifted to matrices in $\Gamma$ via the Euclidean algorithm:
\[
\begin{split}
	& \begin{pmatrix} 0 & -1 \\ 1 & 0 \end{pmatrix},\ldots,\begin{pmatrix} 0 & -1 \\ 1 & 24 \end{pmatrix}, \begin{pmatrix} 1 & 0 \\ 0 & 1 \end{pmatrix}, \\
	& \begin{pmatrix} 1 & 0 \\ 5 & 1 \end{pmatrix}, \begin{pmatrix} -2 & -1 \\ 5 & 2 \end{pmatrix}, \begin{pmatrix} 2 & 1 \\ 5 & 3 \end{pmatrix},
	\begin{pmatrix} -1 & -1 \\ 5 & 4 \end{pmatrix}.
\end{split}
\]
We denote the representatives by $\gamma_1,\ldots,\gamma_{30}$, e.g.\ $\gamma_{26} = I$.

Next we determine the nullspace of the relations matrix. Since the matrix representations of $I+S$, $I+Q+Q^2$ and $I+\varepsilon$ all are $90 \times 90$-matrices, we do not write these up here. After bringing the matrix on echelon form and deleting zero rows, the nullspace matrix of the relations matrix is a $7 \times 90$-matrix.

Next we compute the dimension of the $+$-space. We therefore build an array whose $i$'th entry is the index $j$ of the coset representative $\gamma_j$ satisfying $\gamma_iT = \gamma_0\gamma_j$ for some $\gamma_0 \in \Gamma_0$. This becomes

\begin{center}
\begin{tabular}{ccccccccccccccccc}
	$i$ & \vline & 1 & 2 & 3 & 4 & 5 & 6 & 7 & 8 & 9 & 10 & 11 & 12 & 13 & 14 & 15 \\ \hline
	$j$ & \vline & 2 & 3 & 4 & 5 & 6 & 7 & 8 & 9 & 10 & 11 & 12 & 13 & 14 & 15 & 16 \\ \\
	$i$ & \vline & 16 & 17 & 18 & 19 & 20 & 21 & 22 & 23 & 24 & 25 & 26 & 27 & 28 & 29 & 30 \\ \hline
	$j$ & \vline & 17 & 18 & 19 & 20 & 21 & 22 & 23 & 24 & 25 & 1 & 26 & 27 & 28 & 29 & 30
\end{tabular}
\end{center}

This array shows that $\gamma_1,\ldots,\gamma_{25}$ are in the same cusp equivalence class, and the rest are in their own class. All in all, we see that $\Gamma_0(25)$ has $6$ cusps, represented by $\gamma_1$, $\gamma_{26}$, $\gamma_{27}$, $\gamma_{28}$, $\gamma_{29}$ and $\gamma_{30}$.

We also build an array whose $i$'th entry is the index $j$ of the coset representative $\gamma_j$ satisfying $\gamma_j\varepsilon = \delta_0\gamma_i$ for some $\delta_0 \in \Delta_0(N)$. This becomes:

\begin{center}
\begin{tabular}{ccccccccccccccccc}
	$i$ & \vline & 1 & 2 & 3 & 4 & 5 & 6 & 7 & 8 & 9 & 10 & 11 & 12 & 13 & 14 & 15 \\ \hline
	$j$ & \vline & 1 & 25 & 24 & 23 & 22 & 21 & 20 & 19 & 18 & 17 & 16 & 15 & 14 & 13 & 12 \\ \\
	$i$ & \vline & 16 & 17 & 18 & 19 & 20 & 21 & 22 & 23 & 24 & 25 & 26 & 27 & 28 & 29 & 30 \\ \hline
	$j$ & \vline & 11 & 10 & 9 & 8 & 7 & 6 & 5 & 4 & 3 & 2 & 26 & 30 & 29 & 28 & 27
\end{tabular}
\end{center}

We now take each of the six representatives $\gamma_i$, match them with the corresponding $\gamma_j$ from this table (checking that $\gamma_j$ is in the same equivalence class). We do not need to compute the $\delta_0$'s since we have trivial character, and therefore always will get elements of the $-$-space.

From the table we see that only $\gamma_1$ and $\gamma_{26}$ give rise to $\gamma_j$'s in the same cusp equivalence class, and we therefore get $2$ elements of the $-$-space.

Since there are $6$ cusps, this means that the remaining $4$ are split evenly between the $+$- and $-$-spaces, giving that the dimension of the $+$-space is $2$.

From this we get that the dimension of $S_4(\Gamma_0(25))$ is $7-2 = 5$ (nullspace dimension minus $+$-space dimension).

Finally we compute the Hecke action $t_1,\ldots,t_M$ for $M = \lfloor\frac{30 \cdot 4}{12}\rfloor = 10$.

Each coset representative $\gamma_r$ gives rise to a kernel element $xy \otimes \gamma_r$, and to compute $t_n$ for $xy \otimes \gamma_r$ we compute $A.xy$ for all $A \in H_n$, keeping track of the index $r_A$ of $\gamma_rA = \delta_{0,A}\gamma_{A_r}$.

In the case of $n=3$ we have
\[
	H_3 = \left\{\begin{pmatrix} 1 & 0 \\ 0 & 3 \end{pmatrix}, \begin{pmatrix} 1 & 0 \\ 1 & 3 \end{pmatrix}, \begin{pmatrix} 1 & 0 \\ 2 & 3 \end{pmatrix}, \begin{pmatrix} 2 & 1 \\ 1 & 2 \end{pmatrix}, \begin{pmatrix} 3 & 0 \\ 0 & 1 \end{pmatrix}, \begin{pmatrix} 3 & 1 \\ 0 & 1 \end{pmatrix}, \begin{pmatrix} 3 & 2 \\ 0 & 1 \end{pmatrix}\right\},
\]
and we denote these $A_1,\ldots,A_7$.

We need to compute up to $r=4$ since the forms generated by using $r \leq 3$ only spans a $4$-dimensional space. As a matter of fact, the forms generated by using just $xy \otimes \gamma_4$ span the whole space, and we write a table of indices $r_i$ with respect to $\gamma_4$, as well as the action of $A_i$ on $xy$:

\begin{center}
\begin{tabular}{cccccc}
	$i$ & \vline & 1 & 2 & 3 & 4 \\ \hline
	$r_i$ & \vline & 10 & 22 & 13 & 28 \\ \hline
	$A_i.xy$ & \vline & $3xy$ & $3xy+x^2$ & $3xy+2x^2$ & $2y^2+5xy+2x^2$ \\ \\
	$i$ & \vline & 5 & 6 & 7 \\ \cline{1-5}
	$r_i$ & \vline & 2 & 19 & 11 \\ \cline{1-5}
	$A_i.xy$ & \vline & $3xy$ & $y^2+3xy$ & $2y^2+3xy$
\end{tabular}
\end{center}

Thus we get that
\[
	T_3(xy \otimes \gamma_4) = \sum_{i=1}^7 A_i.xy \otimes \gamma_{r_i},
\]
and we put in the coefficients in a vector $t_3$ (where $(t_3)_{3(r_i-1)+j}$ is the coefficient of $x^jy^{2-j}$ in $A_i.xy$).

We do this for all $n=1,\ldots,10$, and build a matrix with the $t_n$ as columns. We then multiply the nullspace matrix with this matrix, and get (after removing zero rows and putting the matrix on echelon form)
\[
	\begin{pmatrix}
1 & 0 & 0 & 0 & 0 & 0 & 0 & 0 & 1 & 0\\
0 & 1 & 0 & 0 & 0 & 0 & -1 & -1 & 0 & 0\\
0 & 0 & 1 & 0 & 0 & 0 & 1 & -2 & 0 & 0\\
0 & 0 & 0 & 1 & 0 & -1 & 0 & 0 & -3 & 0\\
0 & 0 & 0 & 0 & 1 & 0 & 0 & 0 & 0 & -4
	\end{pmatrix},
\]
which has rank $5$, and the rows therefore form a basis for $S_4(\Gamma_0(25))$. Thus the standard basis (up to $q^{10}$) of this space is
\[
	q+q^9, \ q^2-q^7-q^8, \ q^3+q^7-2q^8, \ q^4-q^6-3q^9, \ q^5-4q^{10}.
\]

\subsection{\texorpdfstring{$S_5(\Gamma_0(12),(\frac{\cdot}{12}))$}{S5(12,(./12))}}
We start out by getting the coset representatives, and we again do this by determining $\P^1(\Z/12\Z)$. We get the $24$ elements
\[
\begin{split}
	& (1,0),\ldots,(1,11),(0,1),(2,1),(2,3),(2,5),(3,1), \\ & (3,2),(3,4),(3,7),(4,1),(4,3),(4,5),(6,1),
\end{split}
\]
and these are lifted to matrices in $\Gamma$ via the Euclidean algorithm:
\[
\begin{split}
	& \begin{pmatrix} 0 & -1 \\ 1 & 0 \end{pmatrix},\ldots,\begin{pmatrix} 0 & -1 \\ 1 & 11 \end{pmatrix}, \begin{pmatrix} 1 & 0 \\ 0 & 1 \end{pmatrix},
	\begin{pmatrix} 1 & 0 \\ 2 & 1 \end{pmatrix}, \begin{pmatrix} 1 & 1 \\ 2 & 3 \end{pmatrix}, \begin{pmatrix} 1 & 2 \\ 2 & 5 \end{pmatrix},
	\begin{pmatrix} 1 & 0 \\ 3 & 1 \end{pmatrix}, \\ & \begin{pmatrix} -1 & -1 \\ 3 & 2 \end{pmatrix}, \begin{pmatrix} 1 & 1 \\ 3 & 4 \end{pmatrix},
	\begin{pmatrix} 1 & 2 \\ 3 & 7 \end{pmatrix}, \begin{pmatrix} 1 & 0 \\ 4 & 1 \end{pmatrix}, \begin{pmatrix} -1 & -1 \\ 4 & 3 \end{pmatrix},
	\begin{pmatrix} 1 & 1 \\ 4 & 5 \end{pmatrix}, \begin{pmatrix} 1 & 0 \\ 6 & 1 \end{pmatrix}.
\end{split}
\]
We denote the representatives by $\gamma_1,\ldots,\gamma_{24}$.

We again determine the relations matrix (of size $184 \times 96$) and its null\-space, which in this case has dimension $8$.

Just as in the last example we choose a coset representative from each cusp equivalence class, and we get $6$ cusps represented by $\gamma_1$, $\gamma_{13}$, $\gamma_{14}$, $\gamma_{17}$, $\gamma_{21}$ and $\gamma_{24}$.

We now take each of these $\gamma_i$, match them with the corresponding $\gamma_j$ satisfying $\gamma_j\varepsilon = \delta_0\gamma_i$ (checking that $\gamma_j$ is in the same equivalence class as $\gamma_i$, which they all are in this case), and compute $\chi(\delta_0)$ to see to which of the $\pm$-spaces they correspond. We summarize the results in the follwing table:

\begin{center}
\renewcommand{\arraystretch}{2}
\begin{tabular}{c|cccccccccccc}
	$i$ & 1 & 13 & 14 & 17 & 21 & 24 \\ \hline
	$j$ & 1 & 13 & 16 & 20 & 23 & 24 \\ \hline
	$\delta_0$ & $\big(\begin{smallmatrix} 1 & 0 \\ 0 & -1 \end{smallmatrix}\big)$ & $\big(\begin{smallmatrix} -1 & 0 \\ 0 & 1 \end{smallmatrix}\big)$ & $\big(\begin{smallmatrix} -5 & 2 \\ -12 & 5 \end{smallmatrix}\big)$ & $\big(\begin{smallmatrix} -7 & 2 \\ -24 & 7 \end{smallmatrix}\big)$ & $\big(\begin{smallmatrix} -5 & 1 \\ -24 & 5 \end{smallmatrix}\big)$ & $\big(\begin{smallmatrix} -1 & 0 \\ -12 & 1 \end{smallmatrix}\big)$ \\ \hline
	$\chi(\delta_0)$ & $-1$ & $1$ & $-1$ & $1$ & $-1$ & $1$ \\ \hline
	Space & $+$ & $-$ & $+$ & $-$ & $+$ & $-$
\end{tabular}
\end{center}

From this we get that the dimension of $S_5(\Gamma_0(12),(\frac{\cdot}{12}))$ is $8-3 = 5$ (nullspace dimension minus $+$-space dimension).

As before we compute the Hecke action $t_1,\ldots,t_M$ for $M = \lfloor\frac{24 \cdot 5}{12}\rfloor = 10$, and in this case it is enough to use the kernel element $(xy^2+xy^2) \otimes \gamma_1$ to generate a basis.







The matrix we get after multiplying the nullspace matrix with the matrix having the $t_n$ as columns (and removing zero rows and putting it on echelon form) is
\[
	\begin{pmatrix}
1 & 0 & 0 & 0 & 0 & 0 & -4 & 0 & -27 & 0\\
0 & 1 & 0 & 0 & 0 & -3 & 0 & -8 & 0 & 0\\
0 & 0 & 1 & 0 & 0 & 0 & -10 & 0 & 12 & 0\\
0 & 0 & 0 & 1 & 0 & -3 & 0 & 0 & 0 & 6\\
0 & 0 & 0 & 0 & 1 & 0 & -5 & 0 & 9 & 0
	\end{pmatrix}.
\]

We thus get that the standard basis (up to $q^{10})$ of $S_5(\Gamma_0(12),(\frac{\cdot}{12}))$ is
\[
	q-4q^7-27q^9, \ q^2-3q^6-8q^8, \ q^3-10q^7+12q^9, \ q^4-3q^6+6q^{10}, \ q^5-5q^7+9q^9.
\]


\begin{thebibliography}{}

\bibitem{Cremona} J. Cremona: \emph{Algorithms for Modular Elliptic Curves}. Cambridge University Press, 1997.

\bibitem{Hardy-Wright} G. Hardy \& E. Wright: \emph{An Introduction to the Theory of Numbers}. 5th ed., Oxford Science Publications, 1979.

\bibitem{Landau} E. Landau: \emph{Elementary Number Theory}. 2nd ed., New York, 1999.

\bibitem{Lucas} E. Lucas: \emph{Sur les Congruences des Nombres Eul\'{e}riens et les Coefficients Diff\'{e}rentiels des Functions Trigonom\'{e}triques Suivant un Module Premier}. Bull. Soc. Math. France 6 (1878), 49--54.

\bibitem{Merel} L. Merel: \emph{Universal Fourier Expansions of Modular Forms} in `On Artin's Conjecture for Odd 2-dimensional Representations'. Lecture Notes in Math. 1585 (1994), Springer, 59--94.

\bibitem{Stanley} R.P. Stanley: \emph{Enumerative Combinatorics, Vol. 1}. Cambridge Studies in Advanced Math. 49, Cambridge Univ. Press, 1997.

\bibitem{Stein} W.A. Stein: \emph{Modular Forms: A Computational Approach}. AMS Graduate Studies in Math. 79, 2007.

\bibitem{Wang} X. Wang: \emph{The Hecke Operators on the Cusp Forms of $\Gamma_0(N)$ with Nebentype} in `On Artin's Conjecture for Odd 2-dimensional Representations'. Lecture Notes in Math. 1585 (1994), Springer, 95--108.

\end{thebibliography}
\end{document}